\documentclass[11pt,leqno]{article}

\usepackage[dvipsnames]{xcolor}
\usepackage{setspace,cite}

\usepackage{float}

\usepackage{amsthm,amsfonts,amssymb,amsmath,oldgerm}
\allowdisplaybreaks

\numberwithin{equation}{section}

\usepackage{enumitem}

\usepackage[utf8]{inputenc}
\usepackage{amssymb} 
\usepackage{amstext}
\usepackage{indentfirst}
\usepackage{amsthm}
\usepackage{fullpage}
\usepackage{oldgerm}
\usepackage{amsmath}
\usepackage{hyperref}
\usepackage{cleveref}
\usepackage{blindtext} \usepackage{titlesec}

\usepackage{mathrsfs}
\numberwithin{equation}{section}

\usepackage{commath} 
\usepackage[mathscr]{euscript} 
\usepackage{comment} 
\usepackage{mathtools} 
\usepackage{verbatim} 
\usepackage{physics} 

\usepackage{wrapfig}
\usepackage{mwe}
\usepackage{epsfig}
\usepackage{color}
\usepackage{graphicx}
\newcommand\R{{\mathbb{R}}} 

\newtheorem{thm}{Theorem}

\newtheorem{prop}{Proposition}
\newtheorem{rem}{Remark}
\newtheorem{definition}{Definition}

\def\phi{\varphi}
\def\k{k}

\newcommand{\RM}{\mathbb{R}}
\newcommand{\ZM}{\mathbb{Z}}

\newcommand{\CM}{\mathbb{C}}

\usepackage{epsfig}
\newcommand{\eq}[2]{\begin{equation}\begin{split}#1\end{split}\label{#2}\end{equation}}

\crefname{lemma}{Lemma}{Lemmas}

\title{Modulational Instability of Small Amplitude Periodic Traveling Waves in the Novikov Equation}

\author{Brett Ehrman\thanks{Department of Mathematics, University of Kansas, 1460 Jayhawk Boulevard, Lawrence, KS 66045, USA; ehrman.brett@ku.edu} ~~ Mathew A. Johnson\thanks{Department of Mathematics, University of Kansas, 1460 Jayhawk Boulevard, Lawrence, KS 66045, USA; matjohn@ku.edu} ~~  St\'ephane Lafortune\thanks{Department of Mathematics College of Charleston Charleston, SC, 29401; lafortunes@cofc.edu}}

\begin{document}
\maketitle

\begin{abstract}
We study the spectral stability of smooth, small-amplitude periodic traveling wave solutions of the Novikov equation, which is a Camassa-Holm type equation
with cubic nonlinearities.  Specifically, we investigate the $L^2(\RM)$-spectrum of the associated linearized operator, which in this case is an integro-differential operator with periodic
coefficients, in a neighborhood of the origin in the spectral plane.   Our analysis shows that such small-amplitude periodic solutions are spectrally unstable to long-wavelength perturbations
if the wave number is greater than a critical value, bearing out the famous Benmajin-Feir instability for the Novikov equation.  
On the other hand, such waves with wave number less than the critical value are shown to be spectrally stable.
Our methods are based on applying spectral perturbation theory to the associated linearization.
\end{abstract}

\section{Introduction}

\noindent
In this paper, our goal is to study the existence and spectral stability of smooth periodic traveling wave solutions in a well-known peakon equation.
The topic of peakon  equations began with the Camassa-Holm (CH) equation
\begin{equation}
u_t-u_{txx}   = 2u_xu_{xx}-3uu_x + uu_{xxx},
\label{CH}
\end{equation}
which was introduced in \cite{Cam,Cam2} as a model of strongly dispersive shallow water waves. Equation \eqref{CH} also appears as a particular case of a family of equations introduced in \cite{Fokas} (see the comment on page 146 of \cite{Fokas95}). Equation \eqref{CH} is integrable in the sense that it is Bi-Hamiltonian, admits infinitely many conserved quantities in involution, has a Lax Pair, 
is solvable by the Inverse Scattering Transform, and passes the Weak Painlev\'e Test \cite{Cam,Cam2,honep,Constantini,Constantini2}. 
Furthermore, (a weak version of) the CH equation admits peakon solutions of the form
\eq{
u=ce^{-|x-ct|},
}{Peakon} 
as well as multi-peakon solutions. 
Historically, the next peakon equation was introduced by Degasperis and Processi \cite{dp} and is given by 
\begin{equation}
u_t-u_{txx}  = 3u_xu_{xx}- 4uu_x + uu_{xxx}.
\label{DP}
\end{equation}
The Degasperis-Processi (DP)  equation was identified in \cite{dp} as being
the only equation, besides the KdV and CH equations, within a family of equations to satisfy asymptotic
integrability conditions up to the third order. Equation \eqref{DP}  was shown  to possess 
the same attributes of integrability as the CH equation, as well as peakon solutions
of the same form \eqref{Peakon},  multi-peakons, and smooth solitons \cite{honep,dhh,Ma05}. A major difference between the DP and the CH equations is that  \eqref{DP}  admits 
discontinuous solutions \cite{Giuseppe1,Giuseppe2,Lundmark1}.
The historically third peakon equation to be introduced  is the following  family parametrized by the parameter $b$ 
\eq{
u_t - u_{xxt} =bu_x u_{xx}-(b+1)uu_x + u u_{xxx},
}{bfamily}  
which is referred to as the b-family.  The $b$-family \eqref{bfamily} was introduced to interpolate between the CH equation
\eqref{CH} ($b=2$) and the DP equation \eqref{DP} ($b=3$) in such a way to preserve the peakon solution \eqref{Peakon} \cite{Holm03,Holm03a}.  
The two integrable cases $b=2,3$ have been singled 
out as the only integrable cases by various tests: the Wahlquist-Estabrook prolongation 
method, the Painlev\'e analysis, symmetry conditions, and a test for asymptotic 
integrability~\cite{dp, wanghone, miknov, honep}. The $b$-family has been derived in the context of the modeling of shallow water waves with $u$ corresponding to the horizontal component of velocity of the fluid \cite{Dullin01,Dullin04,constantin, rossen}(see Proposition 2 in \cite{constantin}, 
and also equation (3.8) in \cite{rossen}).  However, there is a  
debate in the literature about the 
relevancy of the peakon equations in the context of the theory of shallow water waves  \cite{bm}.

One important feature of the three peakon models discussed above is that they all contain quadratic nonlinearities.  
In this paper, we consider the Novikov  equation 
\begin{equation}\label{Novikov}
u_t-u_{xxt}=3uu_xu_{xx} -4u^2u_x+u^2u_{xxx},
\end{equation}
which, historically, appeared as the fourth peakon equation to be considered.
Equation \eqref{Novikov} was proposed by Novikov \cite{Novikov} as part of a classification of polynomial homogeneous generalizations of the Camassa-Holm-type equation with quadratic and cubic nonlinearity that possess an infinite hierarchy of higher-order symmetries. In this regard, \eqref{Novikov} can be regarded as a generalization of the CH equation that accounts for cubic nonlinearities. It is interesting to note however that the Novikov equation only differs from the DP equation \eqref{DP}  by a multiplying factor $u$ applied to the RHS.  Equation \eqref{Novikov} has a Lax pair, is  solvable via the 
Inverse Scattering Transform, has infinitely many symmetries and conserved quantities and is Bi-Hamiltonian \cite{HoneNovikovH,Novikov,HoneAP}.
Further, the Novikov equation has been shown to model the propagation of shallow water waves of moderately large amplitude \cite{Chen2022}.

Similarly to the CH and DP equations and, more generally, to the $b$-family, the Novikov equation admits a variety of both smooth and non-smooth solutions.  
Indeed, in addition to the peakon and multi-peakon solutions, the Novikov equation admits smooth soliton and multi-soliton solutions, as well as smooth and peaked spatially periodic traveling waves
\cite{Matsuno}.   The stability of peaked solutions of \eqref{Novikov} has received considerable attention, 
and their existence, spectral and linear (in)stability, and nonlinear orbital stability have been studied in various works (see \cite{Chen2021b,Chen2019, Liu14,Lafortune2024,Palacio2020,Palacios2021,Wang18}
and references therein).  The orbital stability of smooth soliton has recently been addressed in \cite{Johnson24}. While the stability of smooth periodic solutions to members to the $b$-family has been addressed by several authors \cite{Geyer22,Geyer24,Ehrman24},  to our knowledge, the stability of smooth periodic solutions to the Novikov {equation} has not been studied in the literature.

In terms of well-posedness, the Novikov equation is know to be locally well-posed on $H^s$ for $s>3/2$ in both the periodic and whole line case \cite{Himonas1,Ni2011}, while ill-posed if $s<3/2$ \cite{Himonas2}. Furthermore, in the periodic case, the solutions 
were shown to be global in time for $s>5/2$ if the momentum $m=u-u_{xx}$ does not change sign at $t=0$ \cite{Tiglay} and global for $s>3/2$ if the momentum is nonnegative and $m^{1/3}\in L^2$ \cite{Wu2015}. In \cite{Wu2015}, the global existence of $H^1$ weak periodic solutions is also shown for nonnegative momentum. On the whole line, the weak solutions on $H^1$ were shown to be global if the momentum is initially nonnegative \cite{Wu2011} and strong solutions were shown to be global on $H^s$, $s>3/2$, if the momentum does not change sign at $t=0$ \cite{Wu2012}.

\

The goal of this work is to study the existence of smooth periodic traveling waves of the Novikov equation \eqref{Novikov}, as well as to determine when said waves are 
modulationally stable or modulationally unstable\footnote{By modulational stability / instability, we mean spectral stability / instability of the underlying wave to perturbations in $L^2(\RM)$ in 
a sufficiently small ball about the origin in the spectral plane.  For a precise definition, see Definition \ref{D1} in Section \ref{3} below}.
Modulational stability, as well as spectral stability, has been studied for many dispersive models. Typically, there is a critical number $\k^* > 0$ for which periodic traveling waves of period $2\pi/\k$ are modulationally stable when $\k \in (0,\k^*)$ and modulationally unstable when $\k > \k^*$. {Typically, modulational stability results omit the case where $\k = \k^*$ because the standard approach breaks down in that scenario.

However, a critical wave number does not always exist.} In \cite{J1}, the author found modulational instability criteria for the fractional KdV (Korteweg–De Vries)-type equations that depended exclusively on the parameters in the equation. In particular, the author found that periodic traveling waves with small amplitude of the KdV equation are all modulationally stable. While this is unfortunate in that this means KdV cannot predict the Benjamin-Feir instability of a Stokes wave, it does indicate that such an instability requires either higher-order
dispersive or nonlinear effects to be taken into account. This also means when undertaking modulational instability studies that we should 
be prepared for the possibility that a critical wave number will not exist.
Nevertheless, modulational instability of Stokes waves has been observed in a variety of other nonlinear dispersive models, such
as the generalized Benjamin-Bona-Mahony equation \cite{MH1}, the Camassa-Holm equation\cite{HP1}, and the so-called Whitham equation \cite{HJ15_1,HJ15_2}.

{
\begin{rem}
In addition to the Stokes wave studies described above, there have been a number of works in recent years studying the (spectral) modulational stability and instability
of large amplitude periodic traveling waves in nonlinear dispersive equations: see, for example, \cite{BHJ,JO1,JO2,Conduit}.  Such studies rely on several of the main ideas used in this paper, including
using spectral perturbation theory to encode the critical spectral curves near the origin as roots of a low-dimensional matrix (whose size is determined by the
algebraic multiplicity of the co-periodic generalized kernel of the associated linearized operator).  However, as explicit solutions are rarely known for such large amplitude
waves, the obtained results are typically considerably less explicit than those obtained for Stokes waves (where waves are essentially explicit to {arbitrarily} high order).
Nevertheless, such studies have illuminated deep connections between modulation of periodic waves, their spectral stability, and Whitham's theory of modulations.
\end{rem}
}

In this work, we investigate the modulational instability of such small amplitude periodic traveling waves in the Novikov equation.  Our main result is as follows.

\begin{thm}[Modulational Instability vs. Spectral Stability]\label{T:main}
For any $k>0$, a sufficiently small $2\pi/k$-periodic traveling wave solution of the Novikov equation \eqref{Novikov}, as constructed in Proposition \ref{P:soln_expand},
is modulationally unstable provided that $k^2>3$,  while it is spectrally stable if $k^2<3$.
\end{thm}

More precisely, our analysis in Section \ref{3} yields the following results.  We begin by showing that the limiting equilibrium solution is spectrally stable, meaning all of its
spectrum lies on the imaginary axis.  Specifically, the spectrum of this constant state can be written as a one-parameter family of eigenvalues, all of which are completely confined to
the imaginary axis.  As we bifurcate from the equilibrium solution to nearby small amplitude, $2\pi/k$-periodic traveling waves we utilize a general feature of Hamiltonian stability problems,
which is that a parameter-dependent set of eigenvalues can only bifurcate from the imaginary axis if (at least) two such eigenvalues collide\footnote{To determine whether such a collision
can lead to eigenvalues leaving the imaginary axis, hence instability, one must further investigate the so-called Krein signature of the colliding eigenvalues.  See,
for example, \cite{DT17}}.  As a first preliminary result, we show that if $k^2<3$ then the eigenvalues for the limiting equilibrium state can only collide at the origin in the spectral plane,
and hence for such waves the only possible instability is modulational: see the discussion directly below \eqref{e:ev_cst}.

Continuing, we then show that, near the origin in the spectral plane, the spectrum of the linearization about a sufficiently small $2\pi/k$-periodic traveling wave solution
of the Novikov equation \eqref{Novikov} consists of three spectral curves $\lambda_j(\xi)$ that expand as
\begin{equation}\label{e:spec_curves}
\lambda_j(\xi)=i\xi X_j,~~j=1,2,3
\end{equation}
for $|\xi|\ll 1$, where the $X_j$ are the roots of a cubic polynomial $Q(X)$ with real coefficients that depend 
smoothly on the frequency $k>0$ and the (Bloch) parameter $\xi$: see equation \eqref{e:Q} in Section \ref{3} below.
Our proof shows that for ``long waves," quantified here by $0<k^2<3$, the polynomial $Q(X)$ has three real roots and hence the spectral curves in  \eqref{e:spec_curves} 
are confined to the imaginary axis.  It thus follows that such waves are modulationally stable, and hence spectrally stable by the above abstract Hamiltonian
considerations.  
On the other hand, for ``short waves," quantified by $k^2>3$, the polynomial $Q(X)$ has one real root
and two complex-conjugate roots.  The two complex conjugate roots yield from \eqref{e:spec_curves} the existence of two spectral curves near $\lambda=0$
which bifurcate off of the imaginary axis, implying modulational instability (and hence spectral instability) of the underlying wave.

\

This paper is structured as follows.
In Section \ref{2}, we start by proving existence of smooth periodic traveling wave solutions of equation \eqref{Novikov}. 
We then rescale our coordinate frame to make our solution $2\pi$-periodic and provide an analytic expansion of the solution in terms of its amplitude $a$, culminating in Proposition
\ref{P:soln_expand}.  Our main stability analysis is contained in Section \ref{3}, starting in Section \ref{3.1} with some preliminary results regarding the linearized operators and their Bloch decompositions,
as well as our precise definitions of modulational instability and stability.  In Section \ref{3.2} we use spectral perturbation theory to project the infinite-dimensional
spectral problem onto the three-dimensional critical eigenspace associated with modulational (in)stability in the small amplitude limit, culminating in our proof of Theorem \ref{T:main}.
%
%
%

\

\noindent
\textbf{Acknowledgements:} The work of MJ and BE was partially supported by the NSF under grant
DMS-2108749. The work of SL was supported by a Collaboration Grant for Mathematicians from
the Simons Foundation (award \# 420847).   The authors also thank Wesley Perkins for helpful conversations.

\section{Periodic Traveling Waves: Existence \& Analytic Expansion}
\label{2}

\noindent
We begin by describing the set of periodic traveling wave solutions studied in this work.   Note that  a traveling wave solution of \eqref{Novikov} is a solution
of the form $u(x,t) = \phi(x-ct)$, where the profile $\phi(\cdot)$ satisfies the ODE profile equation
\begin{equation}    \label{existenceODE}
  -c\phi' + c\phi''' + 4\phi^2\phi' = 3\phi\phi'\phi'' + \phi^2\phi'''
\end{equation}
which, after some basic rearranging, can be rewritten as
\begin{equation}\label{existenceODE2}
(\phi^2-c)(\phi-\phi'')'+3\phi\phi'(\phi-\phi'')=0.
\end{equation}
By elementary ODE theory, we necessarily have $\phi\in C^\infty(\RM)$ provided that either $\phi^2(x)<c$ or $\phi^2(x)>c$ for all $x\in\RM$.  Throughout this work,
we will consider waves such that
\begin{equation}\label{condition1}
\phi^2(x)<c~~{\rm for~all}~~x\in\RM.
\end{equation}
Under this assumption, we see that multiplying \eqref{existenceODE2} by the integrating factor $(\phi-\phi'')^{-1/3}$ allows us to rewrite the profile equation 
in the conservative form
\[
\frac{d}{dx}\left(\left(\phi-\phi''\right)^{2/3}\left(c-\phi^2\right)\right)=0.
\]
Additionally taking the condition that
\begin{equation}\label{condition2}
\phi-\phi''>0~~{\rm for~all}~~x\in\RM
\end{equation}
the above can be directly integrated to give
\begin{equation}\label{e:quad2}
\phi-\phi''=\frac{b}{(c-\phi^2)^{3/2}}
\end{equation}
where here $b>0$ is a constant of integration.  Multiplying by $\phi'$ and integrating, \eqref{e:quad2} can be reduced to the quadrature representation
\begin{equation}\label{e:quad}
\frac{1}{2}\left(\phi'\right)^2=E+\frac{\phi^2}{2}-\frac{b\phi}{c\sqrt{c-\phi^2}},
\end{equation}
where here $E\in\RM$ is another constant of integration.

\begin{figure}[t]
\begin{center}
\includegraphics[scale=0.8]{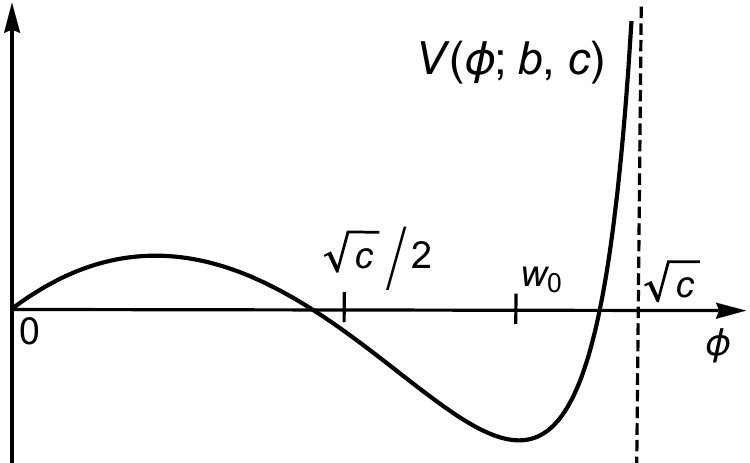}
    \caption{A graphical depiction of the potential $V(\phi;b,c)$ for $c>0$ and $b \in \left(0,\frac{3\sqrt{3}c^2}{16}\right)$.  Phase plane analysis dictates
    the existence of a family of periodic orbits of \eqref{e:quad} surrounding the nonlinear center $w_0$.}\label{GraphV2}
\end{center}    
\end{figure}

By elementary phase plane analysis, the existence of smooth periodic solutions of \eqref{existenceODE} satisfying \eqref{condition1}-\eqref{condition2} follows
provided that the effective potential function
\begin{equation}\label{e:potential}
V(\phi;b,c)=\frac{b\phi}{c\sqrt{c-\phi^2}}-\frac{\phi^2}{2}
\end{equation}
has a strict local minimum in the interval $\phi\in(-\sqrt{c},\sqrt{c})$.  
Conditions guaranteeing the existence of such a minimum
were recently established in \cite[Theorem 2.1]{Johnson24}.  There, it was shown that for each $c>0$ and  $b\in(0,3\sqrt{3}c^2/16)$ the potential satisfies
\[
V(0;b,c)=0,~~V'(0;b,c)>0,~~{\rm and}~~V'(\sqrt{c}/2;b,c)<0
\]
along with
\[
\lim_{\phi\to(\sqrt{c})^-}V(\phi;b,c)=+\infty.
\]
Further, one has 
\[
V''(\phi;b,c)<0~~{\rm for ~all}~~|\phi|<\frac{\sqrt{c}}{2},~~{\rm and}~~V''(\phi;b,c)>0~~{\rm for~ all}~~\frac{\sqrt{c}}{2}<|\phi|<\sqrt{c},
\]
and hence the potential $V$ has a non-degenerate, strict local minimum at some point $w_0\in(\sqrt{c}/2,\sqrt{c})$.  See Figure \ref{GraphV2} for an illustration.
%
By elementary phase plane analysis, it follows that for each $b>0$ and $c>4\cdot 3^{-3/4}b^{1/2}$ there exists a family of periodic solutions of \eqref{e:quad2}
that oscillate about the equilibrium solution $w_0$.

\begin{rem}
The above analysis considers only the case where the profile $\phi$ satisfies both conditions \eqref{condition1} and \eqref{condition2}.  Of course, it is natural 
to consider the existence when these conditions fail.  First, we note that 
the case when $\phi^2<c$ and $m<0$ can be obtained from the case $m>0$ studied above by applying to the solution 
the symmetry $\phi\to-\phi$ of the profile equation \eqref{existenceODE}. Thus the results of this paper extend to the case $m<0$ as well.
Secondly, when $\phi^2>c$ one can show by phase plane analysis
again that there are no smooth periodic traveling wave solutions. Indeed, in that case, the potential is as in \eqref{e:potential}, except that the expression in the radical is $\phi^2-c$. It then can be shown that no value of $b$ can be chosen so that such a function has a local minimum for $\phi>\sqrt{c}$. For further results beyond smooth periodic solutions, see, for example, \cite{PL,ZXO,Li14}.
\end{rem}

Our next goal is to provide an analytic parameterization of the periodic solutions of \eqref{e:quad2} near the 
equilibrium state $w_0$.  To this end, we rescale $z=kx$ for $k>0$ and note that $2\pi/k$-periodic solutions of \eqref{e:quad2} correspond
to $2\pi$-periodic solutions of 
\begin{equation}\label{e:profile4}
-c\phi'+ck^2\phi'''+4\phi^2\phi'=3k^2\phi\phi'\phi''+k^2\phi^2\phi'''
\end{equation}
or, equivalently,
\begin{equation}\label{e:profile3}
\left(c-\phi^2\right)^{3/2}\left(\phi-k^2\phi''\right)=b,
\end{equation}
where now, without loss of generality, the prime notation denotes differentiation in $z$.  
Since the above ODE is invariant under the transformations
$z\mapsto z+z_0$ and $z\mapsto -z$ for all $z_0\in\RM$, there is no loss of generality in seeking solutions of \eqref{e:profile3} that are even in $z$.
With this observation, we let\footnote{Here, $H^2_{\rm per}(0,2\pi)$ denotes the space of $2\pi$-periodic functions in $H^2_{\rm loc}(\RM)$, while 
$H^2_{\rm per,e}(0,2\pi)$ denotes the sector of even functions in $H^2_{\rm per}(0,2\pi)$.}
$F:H^2_{\rm per,e}(0,2\pi)\times\RM^+\times\RM\times\RM\to L^2_{\rm per}(0,2\pi)$ be defined via
\[
F(w;k,b,c):=\left(c-\omega^2\right)^{3/2}\left(\omega-k^2\omega''\right)-b, 
\]
and note that, for $k,b>0$ and $c>4\cdot 3^{-4/3}b^{1/2}$, solutions of \eqref{e:profile3} correspond to solutions $w\in H^2_{\rm per,e}(0,2\pi)$ of
\begin{equation}\label{e:root}
F(w;k,c,b)=0.
\end{equation}
Further, we note that $F$ is an analytic function of its argument

Now, let $k,b>0$ and $c>4\cdot 3^{-4/3}b^{1/2}$ be fixed and note by the above considerations that the local minimum 
$w_0\in(\sqrt{c}/2,\sqrt{c})$ of the effective potential $V(\cdot;b,c)$ in \eqref{e:potential} necessarily satisfies
\[
F(w_0;k,b,c)=0.
\]
By the Implicit Function Theorem, we know that non-constant solutions of \eqref{e:root} may bifurcate
from $w=w_0$ provided that $c>0$ is chosen so that the linearization
\begin{equation}  \label{linearizedODE}
\partial_w F(w_0;k,b,c)=  1-\frac{3w_0^2}{c-w_0^2} -\k^2\partial_z^2
\end{equation}
is not an isomorphism from $H^2_{\rm per}(0,2\pi)$ into $L^2_{\rm per}(0,2\pi)$.
Noting  that
\[
\partial_w F(w_0;k,b,c)\cos(nz)= \left(1-\frac{3w_0^2}{c-w_0^2} +\k^2n^2\right)\cos(nz)
\]
for all $n\in\mathbb{Z}$, it follows that $\cos(z)\in\ker\left(\partial_w F(w_0;k,b,c)\right)$ provided that $c=c_0(b,k)$ where
\[
c_0(b,k)=\left(\frac{k^2+4}{k^2+1}\right)w_0^2.
\]
Further, since the function
\[
\mathbb{N}\cup\{0\}\ni n\mapsto 1-\frac{3w_0^2}{c_0-w_0^2} +\k^2n^2\in\RM
\]
is strictly increasing  in $n$, we see that with the choice of $c_0$ above we actually have
\[
\ker\left(\partial_w F(w_0;k,b,c_0)\right)={\rm span}\left\{\cos(z)\right\}.
\]
Since the equilibrium solution $w_0$ satisfies
\[
\left(c-w_0^2\right)^{3/2}w_0=b
\]
by definition, substituting $c=c_0$ gives the  closed-form expressions\footnote{Observe that $\inf_{k>0}c_0(b,k)=c_0(b,0)=4\cdot 3^{-3/4}b^{1/2}$ and $c_0(b,k)\to\infty$
as $k\to\infty$.  In particular, $c_0(b,k)>4\cdot 3^{-4/3}b^{1/2}$ for all $b>0$ and $k>0$, as needed for the phase-plane based existence theory discussed above.}
\begin{equation}\label{e:equilibrium}
w_0(b,k)=b^{1/4}\left(\frac{3}{k^2+1}\right)^{-3/8}~~{\rm and}~~ ~~c_0(b,k)=b^{1/2}\left(\frac{3}{k^2+1}\right)^{-3/4}\left(\frac{k^2+4}{k^2+1}\right).
\end{equation}

From the above considerations, it follows for all $k>0$ and $b>0$ that the kernel of $\partial_w F(w_0;k,b,c_0)$ is one-dimensional and spanned by precisely
by $\cos(z)$.  Moreover, a straightforward calculation shows that the co-kernel of $\partial_w F(w_0;k,b,c_0)$ is also one-dimensional,
and hence that $\partial_w F(w_0;k,b,c_0)$ is a Fredholm operator with index zero.  Using a straightforward Lyapunov-Schmidt reduction
one may now construct a one-parameter family of non-constant, even, smooth solutions of \eqref{e:root} near the bifurcation
point $(w_0(b,k),c_0(b,k))$.  For details, see \cite{BT,EHJR,J1,kielhofer}, for example.

\begin{prop}\label{P:soln_expand}
For each $k>0$ and $b>0$, there exists a  family of $2\pi/k$-periodic traveling wave solutions of \eqref{Novikov} of the form
\[
u(x,t;a,b)=w\left(k(x-c(k,a)t);a,b,k\right)
\]
for $|a|\ll 1$, where here $w$ and $c$ depend analytically on $a$ and $w$ is a smooth, even and $2\pi$-periodic in $z$, and $c$ is even in $a$.
Furthermore, for $|a|\ll1$ we have the asymptotic expansions
\begin{equation}\label{w V2}
\left\{\begin{aligned}
w(z;a,b,\k) &= w_0(b,\k) + a\cos(z) +a^2\left(d_1+d_2\cos(2z)\right)+ \mathcal{O}(a^3)\\
c(a,b,k)&= c_0(b,\k) + a^2c_2+\mathcal{O}(a^4)
\end{aligned}\right.
\end{equation}
where $w_0$ and $c_0$ are defined above as  in \eqref{e:equilibrium},
\[
d_1=	\frac{\left(1+\k^2\right)^{\frac{5}{8}}\left(5k^4-20\k^2-16\right)}{(48)3^{\frac{5}{8}}b^{\frac{1}{4}}\k^2},~~
d_2=\frac{\left(1+\k^2\right)^{\frac{5}{8}}\left(8+5\k^2\right)}{(12)3^{\frac{5}{8}}b^{\frac{1}{4}}\k^2},
\]
and
\[
c_2=\frac{5}{72}\left(\k^2+4\right)^2.
\]
\end{prop}

Before continuing, it is worth noting that the constant $b>0$ can be scaled to $b=1$ by taking the transformation
\[
w\mapsto b^{1/4}w,~~c\mapsto b^{1/2}c
\]
in the profile equation \eqref{e:quad2}.  We choose not to perform this scaling here, however, as our analysis requires explicit information about the 
variation of both the wavespeed and the profile $w$ in the parameter $b$.  Of course, one could perform the above scaling if desired and determine
these variations through the above transformation. We simply choose to leave the dependence on $b$ explicit throughout.

\section{Spectral Analysis}
\label{3}

\noindent
In this section, we study the stability of the periodic traveling wave solutions constructed in Proposition \ref{P:soln_expand} to small perturbations.
We are specifically interested in the stability {of} so-called localized perturbations, that is, to perturbations in $L^2(\RM)$.  As the boundary 
conditions for the underlying periodic waves and the perturbations are non-commensurate, the stability analysis is rather delicate and
we rely on Floquet-Bloch theory to first characterize the spectrum of the associated linearization in terms of the point spectrum
for a one-parameter family of operators (known as the Bloch operators).  We will then focus our attention on the modulational stability of
the underlying wave, corresponding to the essential spectrum in a sufficiently small neighborhood of the origin in the spectral plane.
By utilizing spectral perturbation theory, we determine sufficient conditions for the modulational instability of the underlying wave.

\subsection{Linearization \& Characterization of the Spectrum}\label{3.1}

\noindent
Throughout this section, let $w=w(\cdot;a,b,k)$ and $c=c(a,b,k)$ for $b,k>0$ and $|a|\ll 1$ form a small amplitude $2\pi/k$-periodic traveling wave solution
of \eqref{Novikov}.  Linearizing \eqref{Novikov} about $w$ in the traveling coordinate frame $z = \k(x-ct)$ yields the linear evolution equation
\begin{equation}\label{e:linear_evol}
    V_t =k\left(1-\k^2\partial_z^2\right)^{-1} \mathcal{L}[w]V
\end{equation}
governing the perturbation $V(z,t)$, where here
\begin{equation}    \label{L1}
\begin{aligned}
    \mathcal{L}[w] &:= c\partial_z - c\k^2\partial_z^3 - 8ww_z - 4w^2\partial_z \\
    &\qquad+ 3\k^2(w_zw_{zz}+ww_{zz}\partial_z+ww_z\partial_z^2) + \k^2\left(w^2\partial_z^3+2ww_{zzz}\right).
\end{aligned}
\end{equation}
Since \eqref{e:linear_evol} is autonomous in time, seeking solutions of the form $V(z,t)=e^{\lambda t}v(z)$
with $\lambda\in\CM$ and $v\in L^2(\RM)$ leads to the spectral problem
\begin{equation}\label{e:spec}
\mathcal{A}[w]v=\lambda v
\end{equation}
where here
\[
\mathcal{A}[w]:=k\left(1-k^2\partial_z^2\right)^{-1}\mathcal{L}[w]
\]
is considered as a closed, densely defined linear operator on $L^2(\RM)$.

We say that the solution $w$ is \emph{spectrally unstable} if the $L^2(\RM)$-spectrum of $\mathcal{A}[w]$ non-trivially
intersects the open right half plane of $\CM$, and is \emph{spectrally stable} otherwise.  Note, however, that since \eqref{e:spec} 
is invariant with respect to the transformations
\[
v\mapsto \overline{v}\quad{\rm and}\quad \lambda\mapsto\overline{\lambda}
\]
and
\[
z\mapsto -z\quad{\rm and}\quad \lambda\mapsto -\lambda,
\]
it follows that the spectrum of $\mathcal{A}[w]$ is symmetric about the real and imaginary axes.  It follows that $w$ is spectrally stable
if and only if 
\[
\sigma\left(\mathcal{A}[w]\right)\subset \RM i
\]
and spectrally unstable otherwise.

As the coefficients of $\mathcal{A}$ are $2\pi$-periodic, it is well-known that non-trivial solutions of \eqref{e:spec} cannot be integrable on $\RM$ and,
at best, that they  can be bounded functions on the line: see, for example, \cite{BHJ,J1}.  Further,  any bounded solution of \eqref{e:spec} is necessarily of the form
\[
v(z)=e^{i\xi z}\phi(z)
\]
where here $\xi\in[-1/2,1/2)$ and $\phi$ is a $2\pi$-periodic function.  From these observations, it  can be shown that $\lambda\in\CM$
belongs to the $L^2(\RM)$-spectrum of $\mathcal{A}[w]$ if and only if there exists a $\xi\in[-1/2,1/2)$ such that the problem
\begin{equation}\label{e:bvp}
\left\{\begin{aligned}
& \mathcal{A}[w]v=\lambda v\\
& v(z+2\pi)=e^{2\pi i\xi}v(z)
\end{aligned}\right.
\end{equation}
has a non-trivial solution, or, equivalently, if and only if there exists a $\xi\in[-1/2,1/2)$ and a non-trivial $\phi\in L^2_{\rm per}(0,2\pi)$ such that
\[
\lambda\phi = e^{-i\xi z}\mathcal{A}[w]e^{i\xi z}\phi =: \mathcal{A}_\xi[w]\phi.
\]  
Furthermore, we have  the spectral decomposition
\begin{equation}
    \sigma_{L^2(\RM)}\left(A[w]\right) = \bigcup_{\xi \in \left(-\frac{1}{2},\frac{1}{2}\right]} \sigma_{L^2_{\rm per}(0,2\pi)}\left(A_{\xi}[w]\right).
    \label{bloch}
\end{equation}
The one-parameter family of operators $\left\{\mathcal{A}_\xi[w]\right\}_{\xi\in[-1/2,1/2)}$ are called the Bloch operators associated to $\mathcal{A}[w]$, and $\xi$ is referred to as
the Bloch frequency.  Since the Bloch operators have compactly embedded domains in $L^2_{\rm per}(0,2\pi)$, it follows for each $\xi\in[-1/2,1/2)$ that the $L^2_{\rm per}(0,2\pi)$
spectrum of $\mathcal{A}_\xi[w]$ consists entirely of isolated eigenvalues with finite algebraic multiplicities.  As a result, the spectral decomposition
\eqref{bloch} provides a continuous parameterization of the essential $L^2(\RM)$-spectrum of $\mathcal{A}[w]$ by a one-parameter family of $2\pi$-periodic eigenvalues
associated to the Bloch operators.	

From above, it follows that to show an underlying wave $w$ is a spectrally stable solution to the Novikov equation \eqref{Novikov} one must determine the eigenvalues
of the Bloch operators $\mathcal{A}_\xi[w]$ for each $\xi\in[-1/2,1/2)$.  Clearly, this is a daunting task.  In this work, we will focus on a particular class of possible
instabilities.

\begin{definition}\label{D1}
A periodic traveling wave solution $w(\cdot;a,b,k)$ of \eqref{Novikov} is said to be modulationally stable if there exists an open neighborhood  
$\mathcal{B}\subset\CM$ of the origin $\lambda=0$ and a $\xi_0>0$ such that
\[
\sigma_{L^2_{\rm per}(0,2\pi)}\left(\mathcal{A}_\xi[w]\right)\cap\mathcal{B}\subset\RM i
\]
for all $|\xi|<\xi_0$.  The wave $w$ is modulationally unstable otherwise.
\end{definition}

In other words, $w$ is modulationally stable provided that it is spectrally stable in a sufficiently small neighborhood of the origin.  As such, it is important to note that 
while modulational instability implies spectral instability, it is not true that modulational stability implies spectral stability since, obviously, there may be additional
instabilities occurring away from the origin.

\begin{rem}
Modulational instabilities are a fundamental feature in many nonlinear systems, including those {arising from} the modeling of surface water waves 
and nonlinear optics.  There has been a considerable amount of work done in connecting such ``spectral" modulational instabilities, as defined above,
to the dynamic instability of periodic traveling wave solutions to slow modulations (via Whitham's theory of modulations).  See, for example,
\cite{BHJ, Conduit,JZ10}.  For more information about Whitham's theory of modulations, see \cite{Whitham} for a physical and mathematical discussion.
\end{rem}

Before moving on, we note that $\phi\in L^2_{\rm per}(0,2\pi)$ solves the problem
\[
\mathcal{A}_\xi[w]\phi=\lambda \phi
\]
if and only if $\phi$ satisfies the problem
\begin{equation}\label{e:g_spec}
k\left(1-k^2(\partial_z+i\xi)^2\right)^{-1}\mathcal{L}_\xi[w]\phi =\lambda \phi,
\end{equation}
where here 
\[
\mathcal{L}_\xi[w]  = e^{-i\xi z}\left(\mathcal{L}[w] \left(e^{i\xi z}~\cdot~\right)\right)
\]
denotes the Bloch operator associated with $\mathcal{L}[w]$ in \eqref{L1}.  
Additionally, we note that the Bloch operators enjoy the spectral symmetries	
\begin{align*}
\sigma\left(\mathcal{A}_\xi[w]\right)=\overline{\sigma\left(\mathcal{A}_\xi[w]\right)}=-\sigma\left(\mathcal{A}_{-\xi}[w]\right)=-\overline{\sigma\left(\mathcal{A}_{-\xi}[w]\right)}.
\end{align*}
It follows that the spectrum of $\mathcal{A}_\xi[w]$ is symmetric about the imaginary axis and, further, that we can without loss of generality restrict ourselves
to positive Bloch frequencies $\xi\in[0,1/2)$.

\subsection{Modulational Instability Analysis}\label{3.2}

\noindent
From the above work, we seek to determine the modulational stability of a given periodic traveling wave $w$ of the Novikov equation
by studying the spectrum of the associated Bloch operators $\mathcal{A}_\xi[w]$ for $|(\lambda,\xi)|\ll 1$.  To this end we begin by first studying
the case $a=0$, in which case the Bloch operators are constant coefficient and hence can be studied directly using Fourier techniques.  We will
then consider the small amplitude solutions $0<|a|\ll 1$ by considering the co-periodic $\xi=0$ case first, followed by a spectral perturbation
argument to study the case $|(\lambda,\xi)|\ll 1$.  

\

To begin, we consider the case $a=0$, corresponding to the constant solution $w=w_0$ in Proposition \ref{P:soln_expand}.  A straightforward Fourier
calculation shows that
\[
\mathcal{A}_\xi[w_0]e^{inz} = \lambda e^{inz}
\]
provided that
\begin{equation}\label{e:ev_cst}
\lambda= \frac{i\left(n+\xi\right)\left(\left(n+\xi\right)^2-1\right)\k^3\left(c_0-w_0^2\right)}{1+\k^2\left(n+\xi\right)^2}=:i\Omega_{n,\xi}
\end{equation}
for some $n\in\ZM$ and $\xi\in[-1/2,1/2)$. Note, in particular, this demonstrates that $w=w_0$ of \eqref{Novikov} is spectrally stable since the $\Omega_{n,\xi}$ are real.
Observe that when $\xi=0$ we have
\[
\Omega_{1,0}=\Omega_{-1,0}=\Omega_{0,0}=0
\]
and, noting that $n\mapsto\Omega_{n,0}$ is odd and strictly increasing in $n$ for $n\geq 2$, we have that
\[
\ldots<\Omega_{-3,0}<\Omega_{-2,0}<0<\Omega_{2,0}<\Omega_{3,0}<\ldots
\]
In particular, it follows that $\lambda=0$ is an isolated eigenvalue of $\mathcal{A}_0[w_0]$ with algebraic multiplicity three 
and
\[
\ker\left(\mathcal{A}_0[w_0]\right)={\rm span}\left\{1,\cos(z),\sin(z)\right\}.
\]
Further, for $\xi\in(0,1/2]$ one can readily verify that when $k^2<3$ we have
\[
\ldots<\Omega_{-3,\xi}<\Omega_{-2,\xi}<\Omega_{0,\xi}<0<\Omega_{-1,\xi}<\Omega_{1,\xi}<\Omega_{2,\xi}<\Omega_{3,\xi}<\ldots.
\]
so that, in particular, the Bloch eigenvalues associated to the constant state $w_0$ never collide away from $(\lambda,\xi)=(0,0)$ 
when the condition $k^2<3$ is satisfied\footnote{ When $k^2>3$ such collisions may occur.  However,
since we are already able to detect a (modulational) spectral instability in this case, we do not study such collisions here.}.  
As such, when $k^2<3$ the only spectral instability possible comes from the $(\lambda,\xi)=(0,0)$ state.  In particular, if we are able to rule
out a spectral (modulational) instability for $|(\lambda,\xi)|\ll 1$ it follows that the associated wave is necessarily spectrally stable.

Our next goal is to track how the eigenvalues $i\Omega_{\pm 1,0}$ and $i\Omega_{0,0}$ bifurcate for $|(a,\xi)|\ll 1$.
To this end, note for $|(a,\xi)|\ll 1$ that the operator $\mathcal{A}_\xi[w(\cdot;a,b,k)]$ is a perturbation of the constant-coefficient operator $\mathcal{A}_0[w_0(b,k)]$ with
\[
\left\|\mathcal{A}_\xi[w(\cdot;a,b,k)-\mathcal{A}_0[w_0(b,k)\right\|_{L^2_{\rm per}(0,2\pi)\to L^2_{\rm per}(0,2\pi)} =\mathcal{O}(|a|+|\xi|)
\]
uniformly in the operator norm.  Consequently, for $|(a,\xi)|\ll 1$ the operator $\mathcal{A}_\xi[w]$ will have three eigenvalues $\lambda_j(a,\xi)=i\Omega_{j,a,\xi}$ for $j=-1,0,1$
that are continuous in $(a,\xi)$ for $|(a,\xi)|\ll 1$ and satisfy $i\Omega_{j,a,\xi}\to 0$ as $(a,\xi)\to (0,0)$.  Further, the three-dimensional total eigenspace
\[
\Sigma_{a,\xi}=\bigoplus_{j\in\{-1,0,1\}}\ker\left(\mathcal{A}_\xi[w]-i\Omega_{j,a,\xi}I\right)
\]
is an analytical continuation of $\Sigma_{0,0}:=\ker\left(\mathcal{A}_0[w_0]\right)$ determined above.  

To track the eigenvalues $\lambda_j(a,\xi)$ for $|(a,\xi)|\ll 1$, our goal is to project the eigenvalue problem \eqref{e:g_spec} onto the three-dimensional
total eigenspace $\Sigma_{a,\xi}$.  More precisely, we aim to construct a suitable basis $\{\phi_j(z;a,\xi)\}_{j=1}^3$ for $\Sigma_{a,\xi}$ and then compute the 
$3\times 3$ matrix
$$
\mathcal{M}_{a,\xi}(\lambda) := \left[\frac{\left<\mathcal{A}_{\xi}[w]\phi_j,\phi_i\right>}{\left<\phi_i,\phi_i\right>} 
	- \lambda \frac{\left<\phi_j,\phi_i\right>}{\left<\phi_i,\phi_i\right>}\right]_{i,j =1,2,3}.
$$
Note that $\mathcal{M}_{a,\xi}(\lambda)$ can be equivalently expressed as
\[
\mathcal{M}_{a,\xi}(\lambda) = \left[\frac{\left<k\mathcal{L}_{\xi}[w]\phi_j,\left(1-k^2(\partial_z+i\xi)^2\right)^{-1}\phi_i\right>}{\left<\phi_i,\phi_i\right>} 
	- \lambda \frac{\left<\phi_j,\phi_i\right>}{\left<\phi_i,\phi_i\right>}\right]_{i,j =1,2,3}
\]
where here $\mathcal{L}_\xi[w]$ and $\mathcal{B}_\xi$ are given as in \eqref{e:g_spec}.  We will find the second formulation above useful
in our forthcoming calculations.
By spectral perturbation theory, it follows that the critical eigenvalues $\lambda_j(a,\xi)$ are then given precisely by the roots of the cubic polynomial
\begin{equation}\label{e:char_poly}
\det\left(\mathcal{M}_{a,\xi}(\lambda)\right)=0.
\end{equation}
It thus remains to find a suitable basis for the total eigenspace $\Sigma_{a,\xi}$, and then to compute the various inner-products above.

To this end, we first note that the above $a=0$ analysis shows that the eigenspace $\Sigma_{0,\xi}$ can be spanned for $|\xi|\ll 1$ by
the  $\xi$-independent orthogonal basis
\[
\phi_1(z)=\cos(z),~~\phi_2(z)=\sin(z),~~\phi_3(z)=1.
\]
We now consider the case $\xi=0$ and $|a|\ll 1$ and aim to construct a basis for the associated total eigenspace
\[
\Sigma_{a,0}={\rm gker}\left(\mathcal{A}_0[w]\right),
\]
which, by the above considerations, must necessarily be 3 dimensional.  
First, owing to the invariance of \eqref{Novikov} with respect to spatial translations we have that
\begin{equation}\label{e:translation}
\mathcal{A}_0[w]w_z=0.
\end{equation}
Alternatively, this can be seen directly by observing that  differentiating the profile equation \eqref{e:profile4}  
with respect to $z$ gives
\[
\mathcal{L}_0[w]w_z=0,
\]
and hence one trivially has \eqref{e:translation}, as claimed.  Similarly, differentiating \eqref{e:profile4} with respect  to the parameters $a$ and $b$ yield
\[
\mathcal{L}_0[w]w_a = -c_a\left(1-k^2\partial_z^2\right)w_z \quad{\rm and}\quad \mathcal{L}_0[w]w_b = -c_b\left(1-k^2\partial_z^2\right)w_z, 
\]
thus yielding
\[
\mathcal{A}_0[w]\left(c_aw_b-c_bw_a\right)=0.
\]
Since the solution $w$ was constructed to be an even function of $z$, by parity it follows that this provides two linearly independent elements
of the kernel of $\mathcal{A}_0[w]$.  Additionally,  by above we have
\[
\mathcal{A}_0[w]w_b= -kc_bw_z,
\]
providing an element of generalized kernel.  Using Proposition \ref{P:soln_expand}, it follows that if we set
\[
d_3 :=2d_1 -\frac{5w_0}{72c_0}\left(\k^2+4\right)^2,
\] 
then the functions
\begin{equation}\label{e:coper_basis}
\left\{
\begin{aligned}
\phi_1(z;a)&:=\frac{2b}{c_0}\left(c_aw_b-c_bw_a\right)=\cos(z)+a\left(d_3+2d_2\cos(2z)\right)+\mathcal{O}(a^2)\\
\phi_2(z;a)&:=-\frac{1}{a}w_z=\sin(z)+2ad_2\sin(2z)+\mathcal{O}(a^2)\\
\phi_3(z;a)&:=\left(\partial_b w_0\right)^{-1}w_b=1+\mathcal{O}(a^2)
\end{aligned}\right.
\end{equation}
provide a normalized basis for the total eigenspace $\Sigma_{a,0}$ for $|a|\ll 1$, and hence also a $\xi$-independent basis
of the critical eigenspace $\Sigma_{a,\xi}$ for $|(a,\xi)|\ll 1$, that continuously extends the basis constructed above for $(a,\xi)=(0,0)$.

To continue, we observe that the operator $\mathcal{L}_\xi[w]$ can be expanded as
\[
\mathcal{L}_\xi[w] = \mathcal{L}_0[w]+i\xi L_1-\frac{1}{2}\xi^2 L_2+\mathcal{O}(\xi)^3
\]
where 
\begin{align*}
\mathcal{L}_0[w] &= (c_0-4w_0^2)\left(\partial_z+\partial_z^3\right)\\
&\qquad+aw_0\left[(2k^2+8)\sin(z)-3k^2\sin(z)\partial_z^2   -(3k^2+8)\cos(z)\partial_z+2k^2\cos(z)\partial_z^3\right]
+\mathcal{O}(a^2)\\
L_1&=(c-4w_0^2)(1+3\partial_z^2)\\
&\qquad+a w_0\left[-\left(3\k^2+8\right)\cos(z)v - 6\k^2\sin(z)\partial_z + 6\k^2\cos(z)\partial_{zz}\right]+\mathcal{O}(a^2)
\end{align*}
and
\[
L_2=-6k^2(c-w_0^2)\partial_z+\mathcal{O}(a^2).
\]
Combining this with the expansions in \eqref{e:coper_basis} yields
\begin{align*}
    \mathcal{L}_{\xi}[w]\phi_1 &= -2\left(c_0-4w_0^2\right)i\xi\cos(z) + a\left[12\k^2\left(c_0-4w_0^2\right)d_2\sin(2z) + w_0\left(10\k^2+16\right)\sin(z)\cos(z)\right]\\
   &+ ia\xi\left[\left(c_0-4w_0^2\right)\left(d_3 - 22d_2\cos(2z)\right)\right]+ ai\xi w_0\left[\left(-9\k^2-8\right)\cos^2(z)+6\k^2\sin^2(z)\right]\\
   & +\xi^2\left[-3\k^2\left(c_0-w_0^2\right)\sin(z)\right] + O\left(a^2 + a\xi^2 + \xi^3\right).
\end{align*}
Additionally, noting that for all $n\in\mathbb{N}$ we have
\begin{align*}
\left(1-k^2(\partial_z+i\xi)^2\right)^{-1}\cos(nz)&=\frac{1+k^2(n^2+\xi^2)}{1+k^4(n^2-\xi^2)^2+2k^2(n^2+\xi^2)}\cos(nz)\\
&\qquad-i\frac{2k^2n\xi}{1+k^4(n^2-\xi^2)^2+2k^2(n^2+\xi^2)}\sin(nz)
\end{align*}
and
\begin{align*}
\left(1-k^2(\partial_z+i\xi)^2\right)^{-1}\sin(nz)&=i\frac{2k^2n\xi}{1+k^4(n^2-\xi^2)^2+2k^2(n^2+\xi^2)}\cos(nz)\\
&\qquad+\frac{1+k^2(n^2+\xi^2)}{1+k^4(n^2-\xi^2)^2+2k^2(n^2+\xi^2)}\sin(nz).
\end{align*}
{{A}} direct calculation using the above asymptotic expansions yield
\begin{align*}
\frac{\left<\mathcal{A}_\xi[w]\phi_1,\phi_1\right>}{\left<\phi_1,\phi_1\right>}&=\left<k\mathcal{L}_\xi[w]\phi_1,\left(1-k^2(\partial_z+i\xi)^2\right)^{-1}\phi_1\right>\\
&=-\frac{2i\xi k (c_0-4w_0^2)}{1+k^2} +O\left(a^2+a\xi^2+\xi^3\right)
\end{align*}
and, similarly,
\begin{align*}
\frac{\left<A_{\xi}[w]\phi_1,\phi_2\right>}{\left<\phi_2,\phi_2\right>} &= k(c_0-4w_0^2)\left(\frac{3}{1+k^2}+\frac{4k^2}{(1+k^2)^2}\right)\xi^2+O\left(a^2+a\xi^2+\xi^3\right)\\
\frac{\left<A_{\xi}[w]\phi_1,\phi_3\right>}{\left<\phi_3,\phi_3\right>} &= ai\xi\left(kd_3(c_0-4w_0^2)-\frac{kw_0\left(3k^2+8\right)}{2}\right)+O\left(a^2+a\xi^2+\xi^3\right).
\end{align*}
By completely analogous calculations, we likewise find 
\begin{align*}
\mathcal{L}_{\xi}[w]\phi_2&= -2i\xi\left(c_0-4w_0^2\right)\sin(z) + \xi^2\left[3\k^2\left(c_0-w_0^2\right)\cos(z)\right]\\
    &+ a\left(c_0-4w_0^2\right)\left[-12d_2\cos(2z)-22d_2\sin(2z)\right]\\
    &+ aw_0\left[-\left(5\k^2+8\right)\cos^2(z) + \left(5\k^2+8\right) \sin^2(z)\right] + ai\xi w_0\left[-\left(15\k^2+8\right)\sin(z)\cos(z) \right]\\
    &+O\left(a^2 + a\xi^2 + \xi^3\right)
\end{align*}
and
\begin{align*}
    \mathcal{L}_{\xi}[w]\phi_3 &= i\xi\left(c_0-4w_0^2\right) + aw_0\left[\left(2\k^2+8\right)\sin(z)\right]\\
    & - ai\xi w_0\left[\left(3\k^2+8\right)\cos(z)\right] + O\left(a^2 + a\xi^2 + \xi^3\right),
\end{align*}
which in turn yield
\begin{align*}
\frac{\left<A_{\xi}[w]\phi_2,\phi_1\right>}{\left<\phi_1,\phi_1\right>} &= -k(c_0-4w_0^2)\left(\frac{3}{1+k^2}+\frac{4k^2}{(1+k^2)^2}\right)\xi^2 +O\left(a^2+a\xi^2+\xi^3\right) = -\frac{\left<A_{\xi}[w]\phi_1,\phi_2\right>}{\left<\phi_2,\phi_2\right>}\\
\frac{\left<A_{\xi}[w]\phi_2,\phi_2\right>}{\left<\phi_2,\phi_2\right>} &= -2i\xi k\frac{(c_0-4w_0^2)}{1+k^2} +O\left(a^2+a\xi^2+\xi^3\right) = \frac{\left<A_{\xi}[w]\phi_1,\phi_1\right>}{\left<\phi_1,\phi_1\right>}\\
\frac{\left<A_{\xi}[w]\phi_2,\phi_3\right>}{\left<\phi_3,\phi_3\right>} &= 0 +O\left(a^2+a\xi^2+\xi^3\right)
\end{align*}
and
\begin{align*}
\frac{\left<A_{\xi}[w]\phi_3,\phi_1\right>}{\left<\phi_1,\phi_1\right>} &= ai\xi\left[2kd_3(c_0-4w_0^2)-\frac{2w_0k^3\left(2k^2+8\right)}{(1+k^2)^2}-\frac{kw_0\left(3k^2+8\right)}{1+k^2}\right] +O\left(a^2+a\xi^2+\xi^3\right)\\
\frac{\left<A_{\xi}[w]\phi_3,\phi_2\right>}{\left<\phi_2,\phi_2\right>} &= a\left[\frac{kw_0\left(2k^2+8\right)}{1+k^2}\right] +O\left(a^2+a\xi^2+\xi^3\right)\\
\frac{\left<A_{\xi}[w]\phi_3,\phi_3\right>}{\left<\phi_3,\phi_3\right>} &= i\xi k(c_0+4w_0^2) +O\left(a^2+a\xi^2+\xi^3\right).
\end{align*}

Taken together, it follows that
\begin{equation}\label{e:Lmatrix}
\begin{aligned}
\mathcal{M}_{a,\xi}(\lambda)&= i\xi
    \begin{bmatrix}
        -2k\alpha m_1 & 0 & 0\\ 0 & -2k\alpha m_1 & 0 \\ & 0 & 0 & k\alpha
    \end{bmatrix}
    + a\begin{bmatrix}
        0 & 0 & 0\\ 0 & 0 & kw_0m_1\left(2\k^2+8\right) \\ 0 & 0 & 0
    \end{bmatrix}\\
    &\quad+ ai\xi\begin{bmatrix}
        0 & 0 & \gamma_1\\
        0 & 0 & 0\\
        \gamma_2 & 0 & 0
    \end{bmatrix}+ \xi^2 \begin{bmatrix}
        0 & k\alpha\left(-3\xi^2m_1+2\xi^2y_1\right) & 0\\ k\alpha\left(3\xi^2m_1-2\xi^2y_1\right) & 0 & 0\\ 0 & 0 & 0
    \end{bmatrix}\\
    &+O\left(a^2+a\xi^2+\xi^3\right),
\end{aligned}
\end{equation}
where, for the sake of notational convenience, we have set
\[
\alpha = c_0-4w_0^2,~~y_1 = -\frac{2k^2}{\left(k^2+1\right)^2},~~ m_1 = \frac{1}{k^2+1}
\]
and
\[
\gamma_1 = 2k\alpha d_3+ky_1w_0\left(2k^2+8\right)-kw_0m_1\left(3k^2+8\right), ~~
\gamma_2 = kd_3\alpha-\frac{kw_0\left(3k^2+8\right)}{2}.
\]

Now, we recall from \eqref{e:char_poly} that the eigenvalues of $\mathcal{A}_\xi[w]$ near $\lambda=0$ in the asymptotic limit $|(a,\xi)|\ll 1$ are precisely
the roots of the cubic polynomial
\begin{equation}\label{e:char_poly2}
\det\left(\mathcal{M}_{a,\xi}(\lambda)\right)=b_0(a,\xi)+ib_1(a,\xi)\lambda+b_2(a,\xi)\lambda^2+ib_3(a,\xi)\lambda^3
\end{equation}
in the variable $\lambda\in\CM$, where the coefficient functions $b_j$, defined for $|(a,\xi)|\ll 1$, depend smoothly on $a$ and $\xi$.  
Now, observing that the spectrum of $\mathcal{A}_\xi[w]$ is symmetric about the imaginary axis, it immediately follows that the $b_j's$ are real-valued
functions.  Further, noting that
$$
\sigma_{L^2_{\rm per}(0,2\pi)}\left(\mathcal{A}_\xi[w]\right) = \overline{\sigma_{L^2_{\rm per}(0,2\pi)}\left(\mathcal{A}_{-\xi}[w]\right)},~~
\sigma_{L^2_{\rm per}(0,2\pi)}\left(\mathcal{A}_\xi[w(\cdot;a)]\right) = {\sigma_{L^2_{\rm per}(0,2\pi)}\left(\mathcal{A}_{\xi}[w(\cdot;-a)]\right)}
$$
it follows that the functions $b_1$ and $b_3$ are even in $\xi$, while the functions $b_0$ and $b_2$ are odd in $\xi$, while all of the 
$c_j's$ are necessarily even in $a$.  Additionally noting that $\lambda=0$ is necessarily a root of \eqref{e:char_poly2} of multiplicity three when
$\xi=0$ and $|a|\ll 1$, it follows that 
\[
b_j(a,\xi)=d_j(a,\xi)\xi^{3-j},~~~j=0,1,2,3,
\]
where the functions $d_j$ are real-valued and depend smoothly on $a$ and $\xi$ for $|(a,\xi)|\ll 1$.  Setting $\lambda=i\xi X$ it follows that
\begin{equation}\label{e:Q}
\det\left(\mathcal{M}_{a,\xi}(\lambda)\right)=i\xi^3\left(d_3(a,\xi)X^3-d_2(a,\xi)X^2-d_1(a,\xi)X+d_0(a,\xi)\right)=i\xi^3 Q(X;a,\xi).
\end{equation}
The underlying wave is thus modulationally stable if the polynomial $Q$ admits three real roots, while it is modulationally unstable otherwise.  

To determine the reality of the roots of $Q$, it is sufficient to study its discriminant
$$
    \Delta(a,\xi) := 18d_3d_2d_1d_0 + d_2^2d_1^2 + 4d_2^3d_0 + 4d_3d_1^3 - 27d_3^2d_0^2,
$$
which can be directly expanded as
\begin{equation}\label{e:disc_expand}
\Delta(a,\xi) = \Delta(0,\xi)+\Lambda(b,k)a^2 + O\left(a^2\left(a^2+\xi^2\right)\right)
\end{equation}
for $|(a,\xi)|\ll 1$.  The polynomial $Q$ will have three real roots, corresponding to modulational stability, when $\Delta(a,\xi)>0$, while it will have one real and two
complex-conjugate roots, corresponding to modulational instability, when $\Delta(a,\xi)<0$.
Using the asymptotic expansion in \eqref{e:Lmatrix}, one finds through direct calculations that
\[
\Delta(0,\xi) = \frac{12\sqrt{3}b^3k^{18}\left(k^2+3\right)^4\left(7k^2+3\right)^2\xi^2}{\left(k^2+1\right)^{\frac{19}{2}}}
\]
and
\[
\Lambda(b,k) = \left(\frac{4b^{\frac{5}{2}}k^{14}\left(k^2+3\right)^3\left(k^2+4\right)^2\left(7k^2+3\right)}{3^{\frac{3}{4}}\left(k^2+1\right)^{\frac{29}{4}}}\right)\left(3-k^2\right).
\]
For $k>0$ and $k\neq\sqrt{3}$, we can directly see that $\Delta(0,\xi)>0$ and that the sign of $\Lambda(b,k)$ is determined from the 
sign of the quantity $3-k^2$.

With the above observations, we can complete the proof of Theorem \ref{T:main}.  By \eqref{e:disc_expand} we see that 
if $k^2>3$ then for fixed but small $|a|$ we can
choose a $\xi_0>0$ such that the discriminant $\Delta(a,\xi)$ is strictly negative for all  $0<|\xi|<\xi_0$.  For such $(k,a)$ it follows that the polynomial $Q$ in \eqref{e:Q} has one real root
and two complex-conjugate roots given by $X_\pm=\alpha\pm i\beta$ with $\alpha,\beta\in\RM$ and $\beta\neq 0$.  By above, these complex conjugate roots yield roots of the 
characteristic polynomial in \eqref{e:char_poly2} of the form
\[
\lambda_{\pm}(\xi) = i\xi\left(\alpha\pm i\beta\right)= \mp\beta \xi + i\alpha\xi.
\]
This immediately implies the modulational instability of periodic traveling wave solutions of \eqref{Novikov} with $k^2>3$ and $|a|$ sufficiently small.

On the other hand, if $k^2<3$ then the discriminant $\Delta(a,\xi)$ is strictly positive for all $|(a,\xi)|\ll 1$.  In this case, it follows that for $|a|$ sufficiently small 
that the roots of $Q$ in \eqref{e:Q} are real and distinct.  It follows that the roots of the characteristic polynomial in \eqref{e:char_poly2} are of the form
\[
\lambda_j(\xi) = i\xi\alpha_j,~~j=1,2,3
\] 
where the $\alpha_j=\alpha_j(a,k,\xi)$ are real and distinct.  
It follows that the  three associated spectral curves of $\mathcal{A}_\xi[w]$ near $\lambda=0$ are confined for $|\xi|\ll 1$ to the imaginary
axis\footnote{ In particular, providing a triple covering of a small interval on $\RM i$ of $\lambda=0$.}, implying modulational stability of the underlying wave.
Recalling that when $k^2<3$ the only possible spectral instability is modulational, i.e. occurs for $|(\lambda,\xi)|\ll 1$, it follows that such waves are necessarily
spectrally stable: see the discussion directly below \eqref{e:ev_cst}.
%
Taken together, this establishes our main result Theorem \ref{T:main}.  

\section{Data Availability Statement} 

Data sharing is not applicable to this article as no new data were created or analyzed in this study.

\bibliographystyle{unsrt}

\end{document}